\documentclass[a4paper,12pt]{article}
\usepackage[dvips]{epsfig}
\usepackage{amsmath,amssymb,amsbsy, amstext,amscd,amsfonts}
\usepackage{color,enumerate,euscript,graphicx,hyperref,srcltx,tikz,indentfirst}
\usepackage{xcolor}
\usepackage{mathtools}
\usepackage{array}
\usepackage{booktabs}
\usetikzlibrary{shapes.geometric}
\usetikzlibrary{ positioning,  shapes.geometric}
\input amssym.def
\input amssym.tex

\newtheorem{lemma}{Lemma}[section]%
\newtheorem{theorem}[lemma]{Theorem}%
\newtheorem{proposition}[lemma]{Proposition}%
\newtheorem{corollary}[lemma]{Corollary}%
\newtheorem{example}[lemma]{Example}%
\newtheorem{question}[lemma]{Question}%

\def\lg{\langle} \def\rg{\rangle} \def\char{\,\hbox{\rm char}\,}
  
\def\Aut{\hbox{\rm Aut\,}} 
  
\def\ZZ{\mathbb{Z}}
 \def\di{\bigm|} \def\oh{\overline H} 
\def\o{{\rm o}}\def\AGL{\hbox{\rm AGL}} 
  
\def\pf{\noindent{\it Proof.} } \def\Syl{\hbox{\rm Syl}}
\def\qed{\hfill $\Box$} \def\demo{\pf}
 
 \def\AGL{\mathrm{AGL}}

 \def\O{\Omega}

   \def\ox{\overline{X}}
\begin{document}
\begin{center}
{\bf\large The classification of Bidihedral Groups}
\footnote{Corresponding author: Pengchong Zeng\\
This research was supported by National Natural Science Foundation of China (12471332,12571362)}
\end{center}


\begin{center}
Hao Yu, Pengchong Zeng\\
\medskip
{\it {\small School of Mathematics \& Center for Applied Mathematics of Guangxi (Guangxi University) \\
Guangxi University, Nanning 530004, P. R. China.}}
\end{center}

\renewcommand{\thefootnote}{\empty}


\begin{abstract}
A group is called bidihedral if it can be expressed as a product of two dihedral subgroups.
In this paper, a complete classification for all bidihedral groups is given.
\end{abstract}

\section{Introduction}\label{Introduction}

A group $G$ is said to be \emph{properly factorizable} if $G = AB$ for two proper subgroups $A$ and $B$ of $G$, while the expression $G = AB$ is called a \emph{factorization} of $G$ and both $A$ and $B$ are \emph{factors} of $G$.
Furthermore, if \( A \cap B = 1 \), then we say that $G$ has an \emph{exact factorization}.

In 1937, Ore posed the following problem\cite{Ore1937}: For a finite group, describe and classify all of its exact factorizations.
Since this problem was proposed, many scholars have conducted extensive research on it, yet it remains unresolved to this day.

The classical works of Itô, Douglas, Wielandt, Kegal, and others laid the foundation for the solvability and structural properties of exact factorizations, demonstrating the powerful potential of decomposing a complex group into two simpler subgroups.
Especially when the factor groups possess favorable properties (such as being abelian, nilpotent, or cyclic), the properties of the whole group can often be strongly controlled, which provides an important reference for understanding more general classes of groups.
In 1955, Itô proved that if a group $G$ can be expressed as the product of two abelian groups, then $G$ is metabelian \cite{Ito1}.
In 1961, Douglas proved that if a group $G$ can be expressed as the product of two cyclic groups, then $G$ is supersolvable \cite{douglas1961supersolvability}.
Also, Wielandt and Kegal proved that if a group $G$ can be expressed as the product of two nilpotent groups, then $G$ is solvable \cite{wielandt1951produkt,Kegel1961}.

Recent interest in exact products of two finite groups has focused on primitive permutation groups containing a regular subgroup{} and factorizations of finite almost simple groups \cite{burness2021solvable, jones2002cyclic, li2023exact, liebeck1996factorizations, liebeck2010regular, xia2017quasiprimitive}.
Hering, Liebeck and Saxl classified factorizations of exceptional groups of Lie type\cite{hering1987factorizations}, and the landmark work by Liebeck, Praeger, and Saxl classified maximal factorizations of almost simple groups \cite{liebeck1996factorizations}.
Recently, all exact factorizations of the almost simple groups have been classified by Li, Wang and Xia\cite{li2023exact}.

A group is \textit{bicycle} if it is a product of two cyclic subgroups.
One of the longstanding unsolved problems in this field is the classification of finite bicyclic groups.
Besides the aforementioned classical result by Douglas on supersolvability of bicyclic groups,
Huppert established that every bicyclic $p$-group is metacyclic if $p>2$ is an odd prime~\cite{Huppert1967}.
However, it is well known that non-metacyclic bicyclic $2$-groups exist.
The metacyclic $p$-groups were classified by Xu and Newman for odd prime $p$~\cite{NX1988}, and by Xu and Zhang
for $p=2$~\cite{XZ2006}.
The bicyclic $2$-groups were eventually classified by Janko~\cite{Janko2008}.
Some results directly related to bicyclic groups can be found in the context of regular embeddings of complete bipartite graphs into orientable surfaces~\cite{Jones2010}, and in the context of Hopf algebras arising from exact factorizations of bicyclic groups~\cite{ABM2014,ACIM2009}.
To the best of our knowledge it remains a challenging open problem in the general case.

The interest in exact factorizations with a cyclic factor extends naturally to combinatorics.
a \textit{Cayley map} $M$ is an embedding of an (undirected, simple, and connected) Cayley graph on  a finite group $H$ into an oriented surface.
It is well known that every Cayley map admits a vertex-transitive subgroup of automorphisms. When this subgroup is also arc-regular, the map is termed a \textit{regular Cayley map}.
It is established that the automorphism group $X:=\Aut(M)$ of a regular Cayley map $M$ admits an exact factorization $X=H\lg\sigma\rg$~\cite{jajcay2002skew,conder2016cyclic}.
Similarly, exact factorizations with a dihedral factor have important applications in unoriented regular Cayley maps~\cite{KK2006, HKK2024, Yu2024}.
Without going into the details, the automorphism group $X:=\Aut(M)$ of an unoriented regular Cayley map $M$ has an exact factorization $X=H\lg c, d\rg$ \cite{KK2006}, where $\lg c, d\rg$ is a dihedral group.

In the study of exact factorizations with cyclic factors, Lucchini provided an upper bound for the cyclic factor (see Proposition \ref{cyclic-point-stabilizer}). However, no such upper bound has been established for exact factorizations with dihedral factors. In this paper, we resolve this gap by establishing a rigorous upper bound for the dihedral factor.

\begin{theorem}\label{dihedral}
Let $G$ be a transitive permutation group of degree $n > 1$ whose point-stabilizer subgroup is dihedral.
Then $|G| < 4n^2$.
In particular, if $G$ is nilpotent, then $|G| \leq n^2$, which is equivalent to $|G| \geq |H|^2$.
\end{theorem}

\noindent\textbf{Remark.}
The bound $n^2$ is the best upper bound as shown by the group $D_8 D_8 = C_2 \wr C_2^2$ of degree $n = 8$ and order $64$.

\vskip 3mm

A group is \textit{bidihedral} if it is a product of two dihedral subgroups;
further, a bidihedral group is exact if the intersection of these two dihedral subgroups is trivial.
Hu and the first author characterized the exact bidihedral groups with order of 4 times an odd integer \cite{yu2025}.
In this paper, a complete classification of all exact bidihedral group is given.
Throughout this paper, let $\lg a,b\rg = \lg a,b \mid a^n = 1,\,b^2 = 1,\,a^b = a^{-1}\rg \cong D_{2n}$ and $\lg c,d\rg = \lg c,d \mid c^m = 1,\,d^2 = 1,\,c^d = c^{-1}\rg \cong D_{2m}$.
Remind that the group $A=G.H$ means that $G\unlhd A$ and for any element $h\in H$, there exists an integer $\ell$ such that $h^{\ell}\in G$.

\begin{theorem}\label{exact-bidihedral}
Let $X = \lg a,b\rg \lg c,d\rg$ be an exact bidihedral group.
Then one of the following forms holds.
\begin{enumerate}
\item $X=( \lg c,d \rg\lg a \rg) \rtimes \lg b \rg$;
\item $X = (\lg a,b \rg \lg c \rg) \rtimes \lg d \rg$;
\item $X = [( \lg c,d \rg\lg a^2 \rg) \rtimes \lg b \rg].\lg a \rg$ or
\item $X = [(\lg a, b \rg \lg c^4 \rg \rtimes \lg d \rg).\lg c^2 \rg].\lg c \rg$.
\end{enumerate}
\end{theorem}

\noindent\textbf{Remark.}
$\lg c,d \rg\lg a \rg,\lg a,b \rg \lg c \rg,\lg c,d \rg\lg a^2 \rg$ and $\lg a, b \rg \lg c^4 \rg$
are exact product groups of the dihedral group and the cyclic group and their structures are classified in \cite{du2025product}.
In particular, all their defining relations are given in \cite{du2025product} provide that the cyclic group is core free.
Thus, we can apply group extension theory to derive all the defining relations of all bidihedral groups by this Theorem.

\section{Preliminaries}

All groups considered in this paper are finite.
For a group $G$ acting on a set $\Omega$, the stabilizer of a point $\alpha \in \Omega$ is denoted by $G_\alpha$.
The core of a subgroup $H$ in $G$, denoted by $H_G$, is the largest normal subgroup of $G$ contained in $H$.
A subgroup $H$ is \emph{core-free} in $G$ if $H_G = 1$.
The \emph{exponent} of a group $G$, denoted $\exp(G)$, is the smallest positive integer $n$ such that $g^n = e$ for all $g \in G$.

\begin{proposition}{\rm \cite{lucchini1998order}}\label{cyclic-point-stabilizer}
Let $G$ be a transitive permutation group of degree $n > 1$ whose point-stabilizer subgroup $H$ is cyclic.
Then $|H|<n$ and $|G| < n^2$, which implies $\frac{|G|}{|H|} > |H|$.
\end{proposition}

\begin{proposition}{\rm\cite[Theorem 1]{M1974}}\label{solvable}
The finite group $G=AB$ is solvable,
where both $A$ and $B$ are subgroups with cyclic subgroups of index no more than 2.
\end{proposition}

\begin{proposition} {\rm\cite[Satz 1]{Ito1}} \label{mateabel}
Let $G=AB$ be a group, where both $A$ and $B$ are abelian subgroups of $G$. Then
$G$ is meta-abelian, that is, $G'$ is abelian.
\end{proposition}

\begin{proposition}{\rm \cite[Corollary 6]{yu2025}}\label{quansi-regular}
Let $X$ be a quasiprimitive permutation group with a stabilizer $H$ and a regular subgroup $K$.
If $X$ is solvable, then one of the the following holds true:
\begin{enumerate}[\rm(a)]
\item If $K$ is cyclic, then either $(X, H, K)=(S_4,D_6,C_4,)$, or $C_p \cong K \leq X \leq \AGL(1, p)$
for some prime $p$.
\item If $K$ is dihedral and $H$ is either cyclic or dihedral, then
$(X, H, K)=(A4,C_3,D_4)$ or $(S_4,D_6,D_4)$.
\end{enumerate}
\end{proposition}

\section{Permutation groups with the special point-stabilizer}
In this section, let $\Omega$ be a finite set with $n$ elements and $\alpha\in\Omega$.
Let $G$ be a transitive permutation group acting on the finite set $\Omega$ and $G_{\alpha}$ a point-stabilizer.
Let $\lg x\rg\leq G_{\alpha}$ such that $k=|G_{\alpha}:\lg x\rg|$.
Then we get $|G:\langle x\rangle| = kn$ and the group $\langle x\rangle$ is also core free in $G$ as $G_{\alpha}$ is core free in $G$.
Now, $G$ is a transitive permutation group acting on the coset $[G:\langle x\rangle]$ with degree $kn$ and
its point-stabilizer is cyclic.
By Proposition~\ref{cyclic-point-stabilizer}, we obtain that $|G_{\alpha}|<k^2n$ and $|G| < (kn)^2$.
Thus, we get the following result.

\begin{lemma}
$|G_{\alpha}|<k^2n$ and $|G| < (kn)^2$.
\end{lemma}

In particular, suppose that $G_{\alpha}$ is dihedral.
We set $G_{\alpha}=\langle c,d \mid c^m = d^2 = 1 ,\, c^d = c^{-1} \rangle\cong D_{2m}$.
Then $\lg c\rg\leq G_{\alpha}$ and $k=|G_{\alpha}:\lg c\rg|$.
Clearly, we get the following corollary, that is the former of Theorem~\ref{dihedral}.
\begin{corollary}
If $G_{\alpha}$ is dihedral, then $|G| < 4n^2$.
\end{corollary}

To proof of the latter of Theorem~\ref{dihedral}, it suffice to show the following lemma.
\begin{lemma}
If $G$ is nilpotent and $G_{\alpha}$ is dihedral, then $|G_{\alpha}| \leq n$, that is $|G_{\alpha}| \leq \frac{|G|}{|G_{\alpha}|}$.
\end{lemma}

\demo
Now suppose that $G$ is nilpotent and $G_{\alpha}$ is dihedral.
Then $G_{\alpha}$ is also nilpotent and so $G_{\alpha}$ is a 2-group.
If $|G| \leq 64$, computation with GAP~\cite{GAP4} shows that the lemma holds.
So in what follows, we assume that $|G|> 64$.
To show this lemma, we can further assume that $|G_{\alpha}|>8$.
Let $N$ be a minimal normal subgroup of $Z(G)$.
Since $Z(G)\neq1$, we get that $N \unlhd G$ , $N \cong \mathbb{Z}_p$ and $G_{\alpha} \cap N = 1$, where p is a prime.

Consider the quotient group $G/N$. Then we have $G_{\alpha}N/N \cong G_{\alpha}/{G_{\alpha} \cap N} \cong G_{\alpha}$ which implies that $G_{\alpha}N/N$ is a dihedral group of order $2m$.
If $G_{\alpha}N/N$ is core-free in $G/N$, then by the induction hypothesis applied to $G/N$,
we obtain $$|G_{\alpha}|=|G_{\alpha}N/N| \leq \frac{|G/N|}{|G_{\alpha}N/N|}=\frac{|G|}{|N||G_{\alpha}|}\leq \frac{|G|}{|G_{\alpha}|},$$
as desired.
So in what follows, we assume that $G_{\alpha}N/N$ is not core-free in $G/N$.

Set $K/N =(G_{\alpha}N/N)_{G/N}\neq1$, the core of $G_{\alpha}N/N$ in $G/N$.
Then $(G/N) / (K/N) \cong G/K$, and $G_{\alpha}K/K$ is core-free in $G/K$.
Using Dedekind's modular law, we get $K = JN$ with $J = G_{\alpha} \cap K$.
Define $\mho_1 = \langle x^p \mid x \in K \rangle$.
Since $\mho_1 \char K \unlhd G$, it follows that $\mho_1 \unlhd G$.
Since $N \cong \mathbb{Z}_p$, for any $x \in J,\ y \in N$, we have $(xy)^p = x^p \in J$, which implies $\mho_1 \leq J$.
Then we get $\mho_1 = 1$ as $G_{\alpha}$ is core free.
Since $N \leq Z(G)$ and $N \cap J \leq N \cap G_{\alpha} = 1$, we get $K = J \times N$.

We claim that $p=2$. Indeed, otherwise, $p$ is an odd prime. Then $J$ is the Sylow 2-subgroup of $K$, and thus $J \char K \unlhd G$, so $J \unlhd G$.
Since $G_{\alpha}$ is core free in $G$, we obtain that $J=1$, hence $K=N$, which implies $K/N=1$, a contradiction.
Since $\mho_1 = \langle x^2 \mid x \in K \rangle=1$, we get that $y^2=1$ for any element $y \in J$.
Thus, $J$ is either $\ZZ_2$ or $\ZZ_2\times\ZZ_2$.
In fact, $J\cong\ZZ_2$. Otherwise, $J = \langle c_1,d \rangle\cong \ZZ_2\times\ZZ_2$, where $c_1\in G_{\alpha}$.
Since $J = G_{\alpha} \cap K \unlhd G_{\alpha}$, we get that $c_1$ is either $c$ or $c^2$, which implies $|G_{\alpha}|=4$ or 8, a contradiction.
Therefore, $J\cong\ZZ_2$. Then $J = \langle c_1 \rangle \leq \langle c \rangle$, where $c_1^2=1$.
Since $|G_{\alpha}|>8$, we get that $G$ is not abelian, which implies $K<G$.
Consider the quotient group  $G/K$.
Since $J = \langle c_1 \rangle\cong\ZZ_2$ and $|G_{\alpha}|>8$, we know $G_{\alpha}K/K$ is dihedral.
Since $|G/K|<|G|$, by the induction hypothesis applied to $G/K$, we obtain
$$
\frac{|G_{\alpha}|^2}{|J|^2} = |\frac{G_{\alpha}N}{K}|^2 \leq |\frac{G}{K}| = \frac{|G|}{|N||J|},
$$
which implies
$$
|G_{\alpha}|^2 \leq \frac{|G||J|}{|N|} = |G|,
$$
as desired.
\qed

Note that the bound $n^2$ is the best upper bound.
See the following example.

\begin{example}
Let $$G=\lg a,b,c,d\mid a^4=b^2=c^4=d^2=a^ba=c^dc=a^da=1,\, a^c=abc^2d,\,b^c=bc^2,\, b^d=a^2b\rg.$$
Checking by GAP~\cite{GAP4}, $G\cong D_8 D_8\cong\ZZ_2 \wr \ZZ_2^2$ is a transitive permutation group acting on the coset $[G:\lg c,d\rg]$ with degree 8 and the point stabilizer $\lg c,d\rg\cong D_8$.
\end{example}

\section{Global constructions of bidihedral groups}

Remind that $\lg a,b\rg = \lg a,b \mid a^n = 1,\,b^2 = 1,\,a^b = a^{-1}\rg \cong D_{2n}$ and $\lg c,d\rg = \lg c,d \mid c^m = 1,\,d^2 = 1,\,c^d = c^{-1}\rg \cong D_{2m}$.
In this section, let $X = \langle a,b\rangle \langle c,d\rangle$ be an exact bidihedral group.
By Proposition~\ref{solvable}, $X$ is solvable.
Let $M$ be a subgroup of $X$ of the largest order such that $\langle c,d \rangle  \leq M \subseteq \langle a \rangle \langle c,d \rangle$.
The following theorem gives a global constructions of $X$.

\begin{theorem}\label{main}
The triple $(M,M_X,X)$ is characterized in Table~\ref{tab:cases}.
\begin{table}[!htp]
\centering
\caption{The forms of M,~M\textsubscript{X}~and~X/M\textsubscript{X}}
\label{tab:cases}
\begin{tabular}{cccc}
\toprule
Case & M & M\textsubscript{X} & X/M\textsubscript{X} \\
\midrule
1 & $\langle a \rangle\langle c,d \rangle$ & $\langle a \rangle\langle c,d \rangle$ & $\mathbb{Z}_2$ \\
2 & $\langle a^2 \rangle\langle c,d \rangle$ & $\langle a^2 \rangle\langle c \rangle$ & $D_8$ \\
3 & $\langle a^2 \rangle\langle c,d \rangle$ & $\langle a^2 \rangle\langle c^2,d \rangle$ & $D_8$ \\
4 & $\langle a^2 \rangle\langle c,d \rangle$ & $\langle a^2 \rangle\langle c^3 \rangle$ & $S_4$ \\
5 & $\langle a^3 \rangle\langle c,d \rangle$ & $\langle a^3 \rangle\langle c^4 \rangle$ & $C_2 \times S_4$ \\
6 & $\langle a^4 \rangle\langle c,d \rangle$ & $\langle a^4 \rangle\langle c^3\rangle$ & $C_2 \times S_4$ \\
7 & $\langle a^4 \rangle\langle c,d \rangle$ & $\langle a^4 \rangle\langle c^4 \rangle$ & $C_2 \wr C_2^2$ \\
\bottomrule
\end{tabular}
\end{table}
\end{theorem}
\demo
All cases with $|X| \leq 64$ are listed in Table~\ref{tab:smallcase}.

\begin{table}[htbp]
\centering
\caption{All cases of $|X| \leq 64$ up to isomorphism}
\label{tab:smallcase}
\begin{tabular}{cccc}
\toprule
$X$ & $M$ & $M_X$ & $X/M_X$ \\
\midrule
$D_4 D_{12}$, $D_8 D_8$, $D_8 D_4$ & $\langle a^2 \rangle\langle c,d \rangle$ & $\langle a^2 \rangle\langle c \rangle$ & $D_8$ \\
$D_4 D_{16}$, $D_4 D_{12}$, $D_4 D_8$, $D_8 D_8$, $D_8 D_4$, $D_{12} D_4$, $D_{16} D_4$ & $\langle a^2 \rangle\langle c,d \rangle$ & $\langle a^2 \rangle\langle c^2,d \rangle$ & $D_8$ \\
$D_4 D_6$ or $D_4 D_{12}$ & $\langle a^2 \rangle\langle c,d \rangle$ & $\langle a^2 \rangle\langle c^3 \rangle$ & $S_4$ \\
$D_6 D_8$ & $\langle a^3 \rangle\langle c,d \rangle$ & $\langle a^3 \rangle\langle c^4 \rangle$ & $C_2 \times S_4$ \\
$D_8 D_6$ & $\langle a^4 \rangle\langle c,d \rangle$ & $\langle a^4 \rangle\langle c^3 \rangle$ & $C_2 \times S_4$ \\
$D_8 D_8$ & $\langle a^4 \rangle\langle c,d \rangle$ & $\langle a^4 \rangle\langle c^4 \rangle$ & $C_2 \wr C_2^2$ \\
else & $\langle a \rangle\langle c,d \rangle$ & $\langle a \rangle\langle c,d \rangle$ & $\mathbb{Z}_2$ \\
\bottomrule
\end{tabular}
\end{table}

So in what follows, we assume $|X| > 64$ and denote $M = \langle a_1\rangle \langle c,d\rangle$, where $a_1\in\lg a \rg$.
Then
$$\langle a_1\rangle = \bigcap_{k,t \in \mathbb{Z}} \langle a_1\rangle^{a^{k}b^{t}} \leq \bigcap_{k,t \in \mathbb{Z}} M^{a^{k}b^{t}} = \bigcap_{r,s,k,t \in \mathbb{Z}} M^{a^{k}b^{t}c^{r}d^{s}} = M_X.$$
Next, we carry out the proof by the following three steps.
\vskip 3mm
{\it Step 1: Case of $M_X\ne 1$.}
\vskip 3mm
Suppose $M_X \neq 1$.
Consider the quotient group $\overline{X} \coloneqq X/M_X = \langle \overline{a},\overline{b}\rangle \langle \overline{c},\overline{d}\rangle$.
We claim that $\langle \overline{a},\overline{b}\rangle \cap \langle \overline{c},\overline{d}\rangle = 1$.
In fact, for any $\overline{a^i b^j}=\overline{c^k d^t}\in \langle \overline{a},\overline{b}\rangle \cap \langle \overline{c},\overline{d}\rangle$, we get $a^i b^j d^{-t} c^{-k}\in M_X$.
Since $M_X \leq M = \langle a_1 \rangle \langle c,d \rangle$, we obtain that $a^i b^j d^{-t} c^{-k} = a_1^x c^y d^z$ for some integers $x,y,z$, which implies $a_1^{-x} a^i b^j = c^y d^z c^k d^t \in \langle a,b \rangle \cap \langle c,d \rangle = 1$.
We obtain that $a^i b^j = a_1^x \in M_X$, thus $c^k d^t \in M_X$, which implies $\overline{a^i}\overline{b^j} = \overline{c^k}\overline{d^t} = 1$.
This establishes our claim that $\langle \overline{a},\overline{b}\rangle \cap \langle \overline{c},\overline{d}\rangle = 1$.

Next, we claim that $M=\lg a^k\rg\lg c, d\rg$, where $k\in\{1,2,3,4\}$.
Otherwise, we know that $\langle \overline{a},\overline{b} \rangle$ is dihedral and $\o(\overline{a})=k \ge 5$.
If $|\lg c,d\rg:\langle c,d \rangle \cap M_X|\leq 2$, then $X/M_X = \langle \overline{a},\overline{b} \rangle \langle \overline{c}, \overline{d} \rangle$, where $|\langle \overline{c}, \overline{d} \rangle|\leq2$.
Since $\langle \overline{a} \rangle\char\langle \overline{a},\overline{b} \rangle \unlhd \overline{X}$, we get $\langle a \rangle M_X \unlhd X$, which implies $\langle a \rangle \langle c,d \rangle \leq X$, a contradiction.
Then both $\langle \overline{a},\overline{b} \rangle$ and $\langle \overline{c},\overline{d} \rangle$ are dihedral.
Let $M_0/M_X$ be the largest subgroup of $\overline{X}$ containing $\langle \overline{c},\overline{d}\rangle$ and contained in the subset $\langle \overline{a}\rangle \langle \overline{c},\overline{d}\rangle$.
By the induction hypothesis applied to $\overline{X}$, we have $M_0/M_X = \overline{\langle a^j\rangle} \langle \overline{c},\overline{d}\rangle$ with $j \in \{1,2,3,4\}$.
Then $M_0 = \langle a^j \rangle \langle c,d \rangle M_X = \langle a^j \rangle M_X \langle c,d \rangle = \langle a_1,a^j\rangle \langle c,d \rangle\leq X$.
By the maximality of $M$, we conclude that $\langle a^j\rangle\leq\langle a_1,a^j\rangle = \langle a_1\rangle$,  a contradiction.

Using the induction hypothesis on $X/M_X$,
noting that $M/M_X$ is core-free in $X/M_X$,
we get that $X/M_X$ is isomorphic to  $\mathbb{Z}_2$, $D_8$, $S_4$, $C_2 \times S_4$ or $C_2 \wr C_2^2$, and correspondingly,
$\o(\overline{a})=k$, where $k\in\{1, 2, 3, 4\}$, and so $a^k\in M_X$.
Since $M=\lg a^i\rg \lg c, d\rg$ and $M_X=\lg a^i\rg \lg c^r\rg$ or $\lg a^i\rg\lg c,d\rg$, we know that $\lg a^i\rg =\lg a^k\rg$,
which implies that $i\in \{1, 2, 3, 4\}$.
Clearly,  if $X/M_X=\ZZ_2$, then $M_X = M = \langle a \rangle \langle c,d \rangle$;
if $X/M_X=D_8$ and $\langle \overline{c},\overline{d} \rangle = \mathbb{Z}_2$, then $M_X$ is either $\langle a^2 \rangle \langle c \rangle$ or $\langle a^2 \rangle \langle c^2,d \rangle$;
if $X/M_X =S_4$ and $\langle \overline{c},\overline{d} \rangle = D_6$, then $M_X=\lg a^2\rg \lg c^3\rg$;
if $X/M_X =C_2 \times S_4$ and $\langle \overline{c},\overline{d} \rangle = D_8$, then $M_X=\lg a^3\rg \lg c^4\rg $;
if $X/M_X =C_2 \times S_4$ and $\langle \overline{c},\overline{d} \rangle = D_6$, then $M_X=\lg a^4\rg \lg c^3\rg $;
and  if $X/M_X=C_2 \wr C_2^2=D_8D_8$ and  $\langle \overline{c},\overline{d} \rangle = D_8$, then $M_X=\lg a^4\rg \lg c^4\rg$.

\vskip 3mm
{\it Step 2:  Show that if $M_X=1$ then $\lg a,b\rg\in \{D_{2kp}\di k=2,3,4$ and $p$ is a prime $\}$.}
\vskip 3mm

Suppose that $M_X = 1$.
Then $M = \langle c,d\rangle$, as $\langle a_1 \rangle \leq M_X$.
If $\o(a) = 2$, by Theorem~\ref{dihedral}, we obtain $|X| < 64$, a contradiction.
Thus, $\o(a)>2$.
Since $\langle a,b\rangle_X \cap \langle a\rangle \char \langle a,b\rangle_X \unlhd X$, we get $\langle a,b\rangle_X =1$.

we consider the faithful (right multiplication) action of $X$ on the set of right cosets $[X:\langle c,d\rangle]$. Then $\langle a,b\rangle$ is regular.
If $X$ is primitive, by Proposition~\ref{quansi-regular}, we have $X = A_4$ or $S_4$ which contradicts $|X| > 64$.
Therefore, $X$ is imprimitive.
Pick a  maximal subgroup $H$ of $X$ which contains $M$ properly.
Then $H=H\cap X=(H\cap \lg a,b\rg)\lg c,d\rg=\lg a^s, b_1\rg \lg c,d\rg<X$,
for some $b_1\in \lg a,b\rg\setminus \lg a\rg$ and some integer $s$.
Using the same argument as $M$, one has $a^s\in H_X$.  Reset $\ox:=X/H_X$.
Consider the faithful primitive action of $\ox$ on $\O_1:=[\ox:\oh]$,
with a cyclic regular subgroup of $\lg\overline{a}\rg$, where $|\O_1|=s$.
By Proposition~\ref{quansi-regular}, we know that $s$ is a prime $p$ such that $\ox\le \AGL(1,p)= \mathbb{Z}_p \rtimes \mathbb{Z}_{p-1}$.
So $\lg \overline{c},\overline{d}\rg \leq \mathbb{Z}_{p-1}$ is cyclic, which implies $c^2 \in H_X$.
In what follows,  we consider these two cases when  $a^s=1$ or $a^s\ne 1$, separately.

\vskip 3mm
Case (1):  $a^s=1$.
\vskip 3mm
In this case, $\o(a) = p$ and $H = \langle b \rangle \langle c,d \rangle$.
Since $|X|>64$, we get $p\ge3$.
If $p = 3$, by Theorem ~\ref{dihedral}, then $|\langle c,d \rangle| < 24$, computation with GAP~\cite{GAP4} shows that it is consistent with our hypothesis.
Now suppose $p \geq 5$.
Clearly, $\lg \overline{a}\rg\lhd\ox$.
Consider the action of $X$ on the set of blocks of length 2 on $\O=[X:\lg c,d\rg]$,
that is the orbital of $\O$ under $H$, with the kernel $H_X$.
If $H_X=1$, then we get $\ox=X$ and $\lg a\rg\lhd X$, a contradiction.
Therefore, $H_X\neq1$.
Then $H_X\nleqq \lg c,d\rg $ (as $\lg c,d\rg_X=1$) so that $H_X$ interchanges two points $\lg c,d\rg$
and $\lg c,d\rg b$, which implies $|H_X/H_X\cap \lg c,d\rg |=2$.
Since $H_X\cap \lg c\rg $ is cyclic and $H_X\cap \lg c\rg $ fixes setwise each block of length 2,
we get $|H_X\cap \lg c\rg |\leq2$.
Since $c^2 \in H_X$, we get that $o(c) \mid 4$.
By Theorem ~\ref{dihedral}, then $|X| < 256$.
Using GAP~\cite{GAP4}, there exists no such group.

\vskip 3mm
Case (2): $a^s\neq1$.
\vskip 3mm

Now assume that $a^p \neq 1$.
Consider the quotient group $\oh:=H/\langle c,d\rangle_H$.
Let $\overline{H_0}$ be the subgroup of the largest order in $\overline{H}$ such that $\langle \overline{c},\overline{d} \rangle  \leq \overline{H_0} \subseteq \overline{\langle a \rangle} \langle \overline{c},\overline{d} \rangle$.

Since $|H| < |X|$, by induction hypothesis, we have $\overline{H_0} = \overline{\langle a^{pk}\rangle} \langle \overline{c},\overline{d}\rangle$ with $k \in \{ 1,2,3,4\}$,
which implies $\lg a^{pk}\rg\lg c,d\rg\lg c,d\rg_H=\lg c,d\rg\lg a^{pk}\rg\lg c,d\rg_H=\lg c,d\rg \lg a^{pk}\rg$,
giving $\lg a^{pk} \rg \lg c,d\rg \le H\le X$.
Therefore, we get $a^{pk}\in M_X$.
Since $M_X=1$ and $a^p=a^s\ne 1$, we get that $\o(a)=kp$ where $ k\in \{2, 3, 4\}.$
Therefore, only the following  three groups are remaining:
$\lg a,b\rg=D_{2kp}$, where $k\in \{2, 3,4\}$ and $p$ is a prime.

\vskip 3mm
{\it Step 3: Show that $\lg a,b\rg$ cannot be $D_{2kp}$,
where $k\in\{2, 3, 4\}$ and $p$ is a prime, provided $M_X=1$.}
\vskip 3mm

Suppose that $\lg a,b\rg\cong D_{2kp}$,  $k\in\{2, 3, 4\}$,
reminding that $H=\lg a^p, b\rg \lg c,d\rg$, $\lg a,b\rg_X=1$, $a^p,c^2\in H_X$, $M_X=1$,
$|X|>64$, $\O=[X:\lg c,d\rg]$ and
$X$ has blocks of length $2k$ acting on $\O$.
Moreover, $\lg a^p\rg _X=1$ and there exists no nontrivial element $a^j\in H$
such that $\lg a^j \rg \lg c,d\rg \le H$.
If $p=2$ or 3, then there exists no such group by GAP~\cite{GAP4}.
Therefore, $p\ge5.$
In what follows,  we consider these three cases when  $k=2,3$ or $4$, separately.

\vskip 3mm
Case (1):  $k=2$.
\vskip 3mm

Consider the quotient group $\oh:=H/\langle c,d \rangle_H$.
Then $\langle \overline{c},\overline{d} \rangle$ is either dihedral or cyclic.
Suppose that  $\langle \overline{c},\overline{d} \rangle$ is dihedral group.
Then calculation with GAP~\cite{GAP4} shows that $H/\langle c,d \rangle_H = D_4 D_6 \cong S_4$, which implies $\langle c,d \rangle_H = \langle c^3 \rangle$.
Since $c^2 \in H_X$ and $\exp(S_4) = 12$, we obtain that $3\mid\o(c)$ and $\langle c^{24} \rangle \leq \lg x^{12}\mid x\in H_X\rg=:\mho_{12}(H_X) \leq \langle c,d \rangle_H=\langle c^3 \rangle$, which implies that $\mho_{12}(H_X)$ is a cyclic group and $\langle c^{24} \rangle \char \mho_{12}(H_X)$.
As $\langle c^{24} \rangle \char \mho_{12}(H_X) \char H_X \lhd X$ and $\langle c,d \rangle_X = 1$ we obtain that $\o(c) \mid 24$, which implies $m=\o(c)=3,6,12$ or $24$.
Suppose that $m\leq12$.
Combing  $\lg a,b\rg_X=1$ with Theorem~\ref{dihedral}, we get that $X$ is a transitive permutation group acting on $2m \leq 24$ points and  $p < 24$.
Since the transitive permutation groups of degree up to 24 are all given in GAP~\cite{GAP4,HULPKE20051},
there exists no such group through the calculation with GAP~\cite{GAP4}, a contradiction.
Now $\o(c) = 24$.
Since $H_X / H_X \cap \langle c,d \rangle_H \cong H_X \langle c,d \rangle_H / \langle c,d \rangle_H \leq H / \langle c,d \rangle_H \cong S_4$, $\langle c^6 \rangle\leq H_X \cap \langle c,d \rangle_H$ and  $c^2\in H_X$, we obtain that $H_X / H_X \cap \langle c,d \rangle_H = S_4$ or $A_4$, which implies that $H_X = C_4.S_4$ or $C_4. A_4$.
Calculation with GAP~\cite{GAP4} shows that both $\Aut(C_4.S_4)$ and $\Aut(C_4.A_4)$ is a $\{2,3\}$-group.
Since $X / C_X(H_X) \lesssim \Aut(H_X)$ and $\o(a^2) = p \geq 5$, we get $a^2 \in C_X(H_X)$.
Then $\langle a^2 \rangle \char \langle a^2 \rangle \times H_X = \langle a \rangle H_X \unlhd X$,
which implies $\lg a^2\rg\unlhd X$, a contradiction again.
Therefore, $\langle \overline{c},\overline{d} \rangle$ is a cyclic group.

Now $\langle \overline{c},\overline{d} \rangle$ is either 1 or $C_2$.
If $\langle \overline{c},\overline{d} \rangle = 1$, then $\langle c,d \rangle \unlhd H$, which implies $\langle a^p \rangle \langle c,d \rangle \leq H$, a contradiction.
Then $\langle \overline{c},\overline{d} \rangle\cong C_2$, and so $H / \langle c,d \rangle_H = \langle \overline{a^p},\overline{b} \rangle \rtimes \langle \overline{c},\overline{d} \rangle = C_2^2 \rtimes C_2$,
where $\langle \overline{c},\overline{d} \rangle\cong C_2$.
Since there exists no nontrivial element $a^j\in H$ such that $\lg a^j \rg \lg c,d\rg \le H$,
we get $\overline{b}\in\lg\overline{a^p},\overline{c},\overline{d}\rg=H/\lg c,d\rg_H\cong D_8$ and
$\lg \overline{c},\overline{d}\rg\cong C_2$ is $\lg \overline{d}\rg$ or $\lg \overline{c}\rg$.

Suppose $\langle \overline{c},\overline{d} \rangle=\lg \overline{d}\rg\cong C_2$.
Then $\overline{d}$ interchange $\overline{a}$ and $\overline{b}$ (or $\overline{ab}$), which implies
$\langle \overline{a^p b} \rangle \lhd \overline{H}$ or $\langle \overline{b} \rangle \lhd \overline{H}$, without loss of generality, denote by $\langle \overline{b_2} \rangle$.
Then $L \coloneqq \langle b_2 \rangle \langle c,d \rangle \leq X$.
$X$ acts on $[X : L]$ by right multiplication, and since $|H : L| = 2$, we obtain that $H$ has an orbit of length $2$ in $[X : L]$.
Since $H_X \unlhd X$, we obtain that all orbits of $H_X$ have the same length either $1$ or $2$, which implies that  $x^2 \in L_X , \forall x \in H_X$.
Since $|L : \langle c,d \rangle| = 2$, we obtain that $g^2 \in \langle c,d \rangle_X = 1 , \forall g \in L$.
Since $c^2 \in H_X$ and $x^2 \in L_X , \forall x \in H_X$, we obtain that $o(c) \mid 8$.
But there exists no such group through the calculation with GAP~\cite{GAP4}, a contradiction.
With the same argument we can exclude the case $\langle \overline{c},\overline{d} \rangle=\lg \overline{c}\rg\cong C_2$.

\vskip 3mm
Case (2):  $k=3$.
\vskip 3mm
In the quotient group $\oh:=H/\langle c,d \rangle_H$, we get that $\overline{\langle a^p,b \rangle} = D_6$ and $\langle \overline{c},\overline{d} \rangle$ is core-free.
If $\langle \overline{c},\overline{d} \rangle$ is cyclic, then $\overline{H} = \langle \overline{a^p},\overline{b} \rangle \rtimes \langle \overline{c},\overline{d} \rangle = D_6 \rtimes C_2$, which implies $\langle \overline{a^p} \rangle \ char\  \langle \overline{a^p},\overline{b} \rangle \unlhd \overline{H}$. Then $\lg a^p\rg\lg c,d\rg\leq H$, a contradiction.
So $\langle \overline{c},\overline{d} \rangle$ is dihedral.
By Theorem~\ref{dihedral}, we obtain that $|\overline{H}| < 144$.
Checking by GAP~\cite{GAP4}, we get that $\oh=D_6D_8\cong C_2 \times S_4$, which implies that $\langle c,d \rangle_H = \langle c^4 \rangle$ and $\langle \overline{a}^p, \overline{b} \rangle$ is core free.

Since $c^2 \in H_X$ and  $\exp(C_2 \times S_4) = 12$,
we get $\langle c^{24} \rangle \leq \mho_{12}(H_X):=\lg x^{12}\mid x\in H_X\rg \leq \langle c,d \rangle_H = \langle c^4 \rangle$.
Since $\mho_{12}(H_X) \ char\ H_X \lhd X$, we get that $\o(c) \mid 24$.
If $\o(c) \leq 12$, then there exists no such group through the calculation with GAP~\cite{GAP4}, a contradiction.
So $\o(c) = 24$, which implies $|X|=288p$.
Computation with GAP~\cite{GAP4} shows that $\oh=H/\langle c^4 \rangle = D_6 D_8 = C_2 \times S_4$ has $2$ normal subgroup isomorphic to $S_4$ and one equals to $\langle \overline{a^p,b} \rangle \langle \overline{c} \rangle$, denote by $\overline{K}$.
Since $\langle c^8 \rangle \char \langle c^4 \rangle = \langle c,d \rangle_H \unlhd H$, we get that $\langle c^8 \rangle \unlhd H$.
Since $C_K(\langle c^8 \rangle) / \langle c^4 \rangle \unlhd K / \langle c^4 \rangle \cong S_4$ and $\overline{c} \in C_K(\langle c^8 \rangle)/\langle c^4 \rangle$ is an element with order 4, we obtain that $C_K(\langle c^8 \rangle) / \langle c^4 \rangle = K / \langle c^4 \rangle$, which means $c^8 \in Z(K)$.
Since $|X : H_X| = |X : K| = 2p$ and $p \ge 5$ is a prime, we get $\Syl_3(X) = \Syl_3(H_X) = \Syl_3(H) = \Syl_3(K)$.
Set $P \in Syl_3(K)$. Then
$$C_3\cong \lg c^8\rg \leq \bigcap_{\forall x \in K} P^x = \bigcap_{\forall x \in X} P^x = P_X.$$
If $|P_X| = 3$, then $P_X = \langle c^8 \rangle \lhd X$, a contradiction.
So $|P_X| = 9$, which implies that $\{P\}=\Syl_3(X)$.
Then $P = \langle a^p, c^8 \rangle\cong C_3\times C_3\unlhd X$.
Therefore, $\langle a^p \rangle \langle c,d \rangle = P \langle c,d \rangle \leq X$, a contradiction again.

\vskip 3mm
Case (3):  $k=4$.
\vskip 3mm

In the quotient group $\oh:=H/\langle c,d \rangle_H$, we get that $\overline{\langle a^p,b \rangle} = D_8$ and $\langle \overline{c},\overline{d} \rangle$ is core-free.
If $\langle \overline{c},\overline{d} \rangle$ is cyclic, then $\overline{H} = \langle \overline{a^p},\overline{b} \rangle \rtimes \langle \overline{c},\overline{d} \rangle = D_8 \rtimes C_2$,
and so $\langle \overline{a^p} \rangle \char  \langle \overline{a^p},\overline{b} \rangle \unlhd \overline{H}$,
which implies $\lg a^p\rg\lg c,d\rg\leq H$, a contradiction.
Therefore, $\langle \overline{c},\overline{d} \rangle$ is dihedral.
Checking by GAP~\cite{GAP4}, we get that $\oh$ is either $D_8D_8\cong C_2 \wr C_2^2$ or $D_8D_6\cong C_2 \times S_4$
Suppose that $H/\langle c,d \rangle_H = D_8 D_8 = C_2 \wr C_2^2$.
Then $\langle c,d \rangle_H = \langle c^4 \rangle$.
Since $\exp(C_2 \wr C_2^2) = 4$, we get $\langle c^8 \rangle \leq \mho_4(H_X):=\lg x^4\mid x\in H_X\rg \leq \langle c^4 \rangle$.
Since $\mho_4(H_X) \char H_X \unlhd X$ and $\langle c,d \rangle_X = 1$, we get that $\o(c) \mid 8$.
But there exists no such group through the calculation with GAP~\cite{GAP4}, a contradiction.
Suppose that $H/\langle c,d \rangle_H = D_8D_6\cong C_2 \times S_4$.
Since $\oh=C_2 \times S_4$ has 2 normal subgroup isomorphic to $S_4$ and one of them contains $\langle \overline{a^p,b} \rangle$, $\o(\overline{a^p})=4$ and $\overline{a^p} \in \overline{H_X}$,
we get $\overline{b} \in \overline{H_X}$, a contradiction again.
\qed

\section{The group extension sequence}

In this section, let $X = \lg a,b\rg \lg c,d\rg$ be an exact bidihedral group.
We shall give the structure of $X$ for each case in Table~\ref{tab:cases}.
Remind that the group $A=G.H$ means that $G\unlhd A$ and for any element $h\in H$, there exists an integer $\ell$ such that $h^{\ell}\in G$.

Case 1: It is easy to see that in this case $X=(\langle a \rangle \langle c,d \rangle) \rtimes \langle b \rangle$.

Case 2: Since $X/M_X = \langle \overline{a}, \overline{b} \rangle \langle \overline{d} \rangle \cong D_8$, we obtain $\langle \overline{a}, \overline{b} \rangle \unlhd X/M_X$, which implies $X = (\langle a,b \rangle \langle c \rangle) \rtimes \langle d \rangle$.

Case 3: Since $M = \lg a^2 \rg \lg c, d \rg$, $M_X = \lg a^2 \rg \lg c^2, d \rg$ and $X/M_X = \lg \overline{a}, \overline{b} \rg \lg \overline{c} \rg \cong D_8$, we obtain $\lg \overline{a} \rg \neq Z(X/M_X)$; otherwise this would contradict the choice of $M$.
Then $Z(X/M_X)$ is either $\lg \overline{b} \rg$ or $\lg \overline{ab} \rg$.
Choose $b$ again such that $Z(X/M_X)=\lg\overline{b}\rg$.
Then $\lg\overline{b}\rg \lg\overline{c}\rg\leq X/M_X$,
which implies that $X/M_X = (\lg\overline{c}\rg \times \lg\overline{b})\rg \rtimes \lg\overline{a}\rg$,
and so $X = [(\lg a^2 \rg \lg c,d \rg) \rtimes \lg b \rg].\lg a \rg$.

Case 4: Since $X/M_X = \lg \overline{a}, \overline{b} \rg \lg \overline{c}, \overline{d} \rg = D_4 D_6 \cong S_4$, we obtain that  $A_4\cong\lg \overline{a},\overline{b} \rg \lg \overline{c} \rg\unlhd X/M_X\cong S_4$, which implies $X = (\lg a, b \rg \lg c \rg) \rtimes \lg d \rg$.

Case 5: Since $M = \lg a^4 \rg \lg c, d \rg$ and $X/M_X = \lg \overline{a}, \overline{b} \rg \lg \overline{c}, \overline{d} \rg = D_6 D_8 \cong C_2 \times S_4$, computations with GAP~\cite{GAP4} show that $S_4\cong\lg \overline{a},\overline{b} \rg \lg \overline{c} \rg\unlhd X/M_X$.
We obtain that $X = (\lg a, b \rg \lg c \rg) \rtimes \lg d \rg$.

Case 6: Since $M = \lg a^4 \rg \lg c, d \rg$ and $X/M_X = \lg \overline{a}, \overline{b} \rg \lg \overline{c}, \overline{d} \rg = D_8 D_6 \cong C_2 \times S_4$,
computations with GAP~\cite{GAP4} show that $S_4\cong\lg \overline{a},\overline{b} \rg \lg \overline{c} \rg\leq\overline{X}$.
We obtain that $X = (\lg a, b \rg \lg c \rg) \rtimes \lg d \rg$.

Case 7: Since $M = \lg a^4 \rg \lg c, d \rg$ and $X/M_X = \lg \overline{a}, \overline{b} \rg \lg \overline{c}, \overline{d} \rg = D_8 D_8 \cong C_2 \wr C_2^2$, we obtain that $\lg \overline{a} \rg \lg \overline{c}, \overline{d} \rg$ and $\lg \overline{a^2} \rg \lg \overline{c}, \overline{d} \rg$ is not a subgroup of $\overline{X}$, otherwise this would contradict the choice of $M$.
Computation with GAP~\cite{GAP4} shows that
$$\lg \overline{a},\overline{b} \rg \lg \overline{d} \rg\leq \lg \overline{a},\overline{b} \rg \lg \overline{d} \rg \lg \overline{c^2} \rg\leq X/M_X,$$
where $\overline{d_1} \in \{\overline{d}, \overline{cd}, \overline{c^2 d}, \overline{c^3 d}\}$.
Choose $d$ again such that $d_1=d$.
We obtain that $\overline{X} = [(\lg \overline{a},\overline{b} \rg \rtimes \lg \overline{d} \rg) \rtimes \lg \overline{c^2} \rg] . \lg \overline{c} \rg$, which implies $X = [(\lg a, b \rg \lg c^4 \rg \rtimes \lg d \rg) . \lg c^2 \rg] . \lg c \rg$.

To sum up, we get that $X$ is $( \lg c,d \rg\lg a \rg) \rtimes \lg b \rg$,
$(\lg a,b \rg \lg c \rg) \rtimes \lg d \rg$, $X = [( \lg c,d \rg\lg a^2 \rg) \rtimes \lg b \rg].\lg a \rg$ or
$[(\lg a, b \rg \lg c^4 \rg \rtimes \lg d \rg).\lg c^2 \rg].\lg c \rg$.
Then we finish the proof of the Theorem~\ref{exact-bidihedral}.

Note that
$\lg c,d \rg\lg a \rg,\lg a,b \rg \lg c \rg,\lg c,d \rg\lg a^2 \rg$ and $\lg a, b \rg \lg c^4 \rg$
are exact product groups of the dihedral group and the cyclic group and their structures are classified in \cite{du2025product}.
In particular, all their defining relations are given in \cite{du2025product} provide that the cyclic group is core free.
Thus, we can apply group extension theory to derive all the defining relations of all exact bidihedral groups by  Theorem~\ref{exact-bidihedral} and the group extension order is the following theorem.
\begin{theorem}
All group extension orders of $X$ are the following.
\begin{enumerate}
\item $X = (\lg a \rg_{\lg c,d \rg\lg a \rg}.\lg c,d \rg\lg a \rg) \rtimes \lg b \rg$;
\item $X = (\lg c \rg_{\lg a,b \rg \lg c \rg}.\lg a,b \rg \lg c \rg) \rtimes \lg d \rg$;
\item $X = [(\lg a^2 \rg_{\lg c,d \rg\lg a^2 \rg}.\lg c,d \rg\lg a^2 \rg) \rtimes \lg b \rg].\lg a \rg$ or
\item $X = [((\lg c^4 \rg_{\lg a, b \rg \lg c^4 \rg}.\lg a, b \rg \lg c^4 \rg) \rtimes \lg d \rg).\lg c^2 \rg].\lg c \rg$.
\end{enumerate}
\end{theorem}

\section{On bidihedral groups}

In this section, let $X=\lg a,b \rg \lg c,d \rg$ be a bidihedral group.
If $\lg a,b\rg\cap\lg c,d\rg=1$, then $X$ is given in Theorem~\ref{exact-bidihedral}.
So in what follows, we assume that $\lg a,b\rg\cap\lg c,d\rg\neq1$.
Then $\lg a,b\rg\cap\lg c,d\rg$ is either $1<\lg c_0\rg\leq \lg c\rg$ or $\lg c_1,d_1\rg\leq \lg c,d\rg$, where $d_1\in \lg c,d\rg\setminus\lg c\rg$.

Suppose that $\lg a,b\rg\cap\lg c,d\rg=\lg c_0\rg\neq1$.
Then either $\lg c_0\rg\leq \lg a\rg$ or $c_0=b_0$ where $b_0\in\lg a,b\rg\setminus\lg a\rg$.
If $\lg c_0\rg\leq \lg a\rg$, then in $X/\lg c_0\rg=\lg \overline{a},\overline{b} \rg \lg \overline{c}, \overline{d} \rg$, we get $\lg \overline{a},\overline{b} \rg \cap\lg \overline{c}, \overline{d} \rg=1$,
which implies that $X/\lg c_0\rg$ is given in Theorem~\ref{exact-bidihedral} or $X/\lg c_0\rg=\lg\overline{a},\overline{b}\rg\rtimes\lg\overline{d}\rg$.
If $c_0=b_0$, then $X=\lg a\rg\lg c,d\rg=\lg c,d\rg\lg a\rg$ is given in \cite{du2025product}, where $\lg a\rg\cap \lg c,d\rg=1$.

Suppose that $\lg a,b\rg\cap\lg c,d\rg=\lg c_1,d_1\rg\leq \lg c,d\rg$, where $d_1\in \lg c,d\rg\setminus\lg c\rg$.
Then $\lg c,d\rg=\lg c,d_1\rg$ and $X=\lg a,b\rg\lg c\rg$, where $\lg a,b\rg\cap\lg c\rg=\lg c_1\rg$.
If $c_1=1$, then $X=\lg a,b\rg\lg c\rg$ is given in \cite{du2025product}, where $\lg a,b\rg\cap \lg c\rg=1$.
If $1<\lg c_1\rg\leq \lg a\rg$, then similarly we get that $X/\lg c_1\rg$ is given in Theorem~\ref{exact-bidihedral} or $X/\lg c_1\rg=\lg\overline{a},\overline{b}\rg$ is dihedral.
Now suppose that $c_1=b_1$, where $b_1\in\lg a,b\rg\setminus\lg a\rg$.
Then $X=\lg a\rg\lg c\rg$, where $\lg a\rg\cap\lg c\rg=1$ and there exists an element $d_2\in\lg c,d\rg\setminus\lg c\rg$ such that $d_2\in \lg a\rg$.
If $\o(a)=4$, then $X=\lg a\rg\lg c\rg$ and $\lg c\rg_X=\lg c\rg$ or $\lg c^2\rg$,
which implies $X=(\lg c\rg\rtimes\lg a^2\rg).\lg a\rg$.
Now suppose that $\o(a)>4$.
Then $b_1\notin X'$.
By Proposition~\ref{mateabel}, we get that $X'$ is abelian.
Since $\lg a^2,c^2\rg\leq X'$ and $b_1=c_1\in\lg c\rg$,
we get $\o(c^2)$ is odd and $\lg c\rg=\lg b_1\rg\times\lg c^2\rg$.
Then $X=\lg a,b\rg\lg c^2\rg$ is given in \cite{du2025product}, where $\lg a,b\rg\cap \lg c^2\rg=1$.

Now we get the following Theorem.
\begin{theorem}
Let $X=\lg a,b \rg \lg c,d \rg$ be a bidihedral group.
Then one of the following holds.
\begin{enumerate}
\item $X$ is an exact bidihedral group.
\item If $1\neq\lg a,b\rg\cap\lg c,d\rg\leq \lg c\rg$, then one of the following holds.
\begin{enumerate}
\item $X/\lg a,b\rg\cap\lg c,d\rg$ is an exact bidihedral group or an exact product group of the dihedral group and the cyclic group.
\item $X$ is an exact product group of the dihedral group and the cyclic group.
\end{enumerate}
\item If $\lg a,b\rg\cap\lg c,d\rg=\lg c_1,d_1\rg\leq \lg c,d\rg$, then one of the following holds.
\begin{enumerate}
\item $X/\lg c_1\rg$ is either an exact bidihedral group or a dihedral group.
\item $X$ is an exact product group of the dihedral group and the cyclic group.
\item $X=(\lg c\rg\rtimes\lg a^2\rg).\lg a\rg$ with $o(a) = 4$.
\item $X=\lg a,b\rg\lg c^2\rg$ is an exact product group of the dihedral group and the cyclic group.
\end{enumerate}
\end{enumerate}
\end{theorem}
Remind that exact product groups of the dihedral group and the cyclic group are given in \cite{du2025product}.

\section{the Further Research}

In this paper, we classify all exact bidihedral groups.
The exact product group is related to regular unoriented Cayley maps.
In particular, the automorphism group of regular unoriented Cayley maps of dihedral groups is an exact bidihedral groups, which is given in Theorem~\ref{exact-bidihedral}.

\begin{question}
Classify all regular unoriented Cayley maps of dihedral groups.
\end{question}

{\it Email address:}

Hao Yu: haoyu@gxu.edu.cn

Zeng PengChong: 2406301059@st.gxu.edu.cn

\end{document}